\documentclass[11pt]{article}
\usepackage{a4}
\usepackage{graphicx}
\usepackage{parskip}

\usepackage{amsfonts}
\usepackage{natbib}
\usepackage{algorithm}
\usepackage{algorithmic}

\usepackage{latexsym}

\newtheorem{theorem}{Theorem}[section]

\newtheorem{lemma}{Lemma}[section]
\newtheorem{corollary}{Corollary}[section]
\newtheorem{remark}{Remark}[section]

\newtheorem{condition}{Condition}[section]
\newtheorem{definition}{Definition}[section]

\newcommand{\R}{\mathbb{R}}

\newcommand{\PP} {{  \rm I\hskip-0.22em P}}

\newcommand{\EE} {{\rm I\hskip-0.48em E}}

\usepackage{natbib}

\begin{document}

\centerline{\bf Generic chaining and the $\ell_1$-penalty}

\centerline{Sara van de Geer}

{\bf Abstract} We address the choice of the tuning parameter
$\lambda$ 
in $\ell_1$-penalized M-estimation. Our main concern is models which are
highly nonlinear, such as the Gaussian mixture model. The number of parameters $p$
is moreover large, possibly larger than the number of observations $n$.
The generic chaining technique of \cite{talagrand2005generic}
is tailored for this problem.  It leads to the choice $\lambda \asymp \sqrt {\log p / n }$,
as in the standard Lasso procedure (which concerns the linear model and least squares loss).

\section{Introduction}\label{introduction.section}
Let $X_1 , \ldots , X_n$ be independent observations with values
in some observation space ${\cal X}$, and let for $\theta $ in a parameter
space $\Theta \subset \R^p$ be
given a loss function
$\rho_{\theta} : {\cal X } \rightarrow \R $.
The parameter $\theta$ is potentially high-dimensional, i.e.\ possibly
$p \gg n $. 
In this article we study the $\ell_1$-regularized M-estimator
$$\hat \theta := \arg \min_{\theta \in \Theta} \biggl \{ P_n \rho_{\theta} + \lambda \| \theta \|_1 \biggr \} . $$
Here, we use the notation $P_n \rho_{\theta} :=  \sum_{i=1}^n \rho_{\theta} ( X_i) /n$, 
i.e., it is the empirical measure of the loss function $\rho_{\theta}$, often referred to
as the empirical risk. Moreover, $\lambda>0$ is a tuning parameter and
$\| \theta \|_1 := \sum_{j=1}^p |\theta_j | $ is the $\ell_1$-norm of $\theta$. 

A special case is the Lasso (\cite{tibs96}), which has quadratic loss:
$$\rho_{\theta} ( X) := ( Y - \theta^T Z )^2, \ X= (Y,Z) , $$
where $Y \in \R$ is the response variable and $Z\in \R^p$ are covariables.
There are many papers on the Lasso, see for example \cite{vdG:2001}, \cite{Bunea:06},
\cite{Bunea:07a}, \cite{buneaetal06}, 
\cite{vandeG08}, \cite{koltch09a}, \cite{bickel2009sal}. For an overview and further
results, see also \cite{BvdG2011}.
It is known that generally the choice 
$\lambda \asymp \sqrt {\log p / n } $
is appropriate. Under some distributional assumptions.
this choice leads to favorable theoretical properties of the Lasso,
such as good oracle bounds for the estimation and prediction error.

In this note we address the following question: is
the choice $\lambda \asymp \sqrt {\log p / n }$ also appropriate for
non-linear situations? 
 The example described above is a linear situation. More generally, we call the situation linear
 if 
 for some $\psi : {\cal X} \rightarrow \R^p $ 
$$(P_n - P) ( \rho_{\theta} - \rho_{\tilde \theta} ) = ( \theta - \tilde \theta)^T (P_n-P) \psi, \   \forall \ \theta , \ \tilde \theta ,$$
where
$P \rho_{\theta} := {1 \over n} \sum_{i=1}^n \EE \rho_{\theta} (X_i) $ is the theoretical risk. 
Any generalized linear model
(GLM)
loss function with canonical link function and fixed design is a linear situation. Also density estimation using
an exponential family is a linear situation.
A non-linear situation occurs for
instance in linear least squares regression with random design. Our focus is more on 
other examples, such as mixture models (\cite{Staedler:10}) or mixed effect models (\cite{schelldorfer2011estimation}),
etc.

Let us define the ``true" parameter
$$\theta^0 := \arg \min_{\theta \in \bar \Theta} P \rho_{\theta} ,$$
where we assume $\rho_{\theta}$ is defined
for all $\theta$ in the possibly extended space $\bar \Theta \supset \Theta$.
Let $\theta^* \in \Theta$ be some ``approximation" of $\theta^0$. Here, we have in mind
the best approximation within $\Theta$ (in the case of misspecified models),
and possibly the best ``sparse" approximation (see Remark \ref{theta*.remark} for a definition).
Our choice for the tuning parameter is governed by the behavior over
$\ell_1$-balls $\Theta_M ( \theta^*) := \{ \theta \in \Theta_*: \ \| \theta - \theta^* \|_1 \le M \} $ of the empirical process
$(P_n - P) ( \rho_{\theta}- \rho_{\theta_*} ) $, where $\Theta_* = \Theta$ or, in the case of 
Theorem \ref{oracle2.theorem} (convex loss)
$\Theta_*$ is the smallest convex set containing $\Theta$. 

In the linear case, the supremum of the empirical process can be easily bounded using the dual norm
inequality
\begin{equation}\label{dual-norm.equation}
\sup_{\theta \in \Theta_M ( \theta^*)  } | (P_n - P) ( \rho_{\theta} - \rho_{\tilde \theta} ) |
\le \| (P_n-P) \psi \|_{\infty} M , 
\end{equation}
where for a vector $v \in \R^p $, $\| v \|_{\infty} := \max_{1 \le j \le p } | v_j | $ is the uniform norm.
Moreover, for example for ${\cal N} (0,1/n)$-random variables
$\{ V_j \}_{j=1}^p $ (say), it holds that
$$ \max_{1 \le j \le p } | V_j | = {\mathcal O}_{\PP} \biggl (  \sqrt { \log p \over n } \biggr )  . $$

 We show in this paper that in many non-linear cases, one still has
\begin{equation}\label{supremum.equation}
\sup_{\theta \in \Theta_M (\theta^*) } | (P_n - P) ( \rho_{\theta} - \rho_{\theta^*} ) | =
{\mathcal O}_{\PP} \biggl ( \sqrt {\log p \over n } \biggr ) M.  
\end{equation}
This follows rather easily from a generic chaining (\cite {talagrand1996majorizing})
 and Sudakov minoration argument. We will use the book \cite{talagrand2005generic}.

In the case of regression with robust GLM loss (robust quasi-likelihood loss  functions, quantile functions),
we have
$$\rho_{\theta} (X) = \rho ( Y , \theta^T Z) ,$$
with $\rho(y, \cdot) $ Lipschitz for all $y$. In that situation, one may apply
the contraction inequality (\cite{Ledoux:91}) to arrive at (\ref{supremum.equation}). We will explain this in Subsection \ref{contraction.section}.

Our emphasis is however on cases that go beyond GLM loss.
An example is the Gaussian mixture model
$$ \rho_{\theta} (Y,Z) = \log \biggl (  \sum_{k=1}^r \pi_k \phi_{\sigma_k} ( Y- \beta_k^T Z_k ) \biggr ) , $$
where
$\phi_{\sigma} = \phi ( \cdot/ \sigma ) / \sigma$ is the density of the ${\cal N} (0, \sigma^2)$-distribution,
$\{ \pi_k\}_{k=1}^r$ are mixing coefficients ($\sum_{k=1}^r \pi_k = 1 $), 
$\beta_k $ and $Z_k$ are vectors in $\R^{p_k}$, $k=1 , \ldots , r$, and where
$Y$ is again a response variable and 
$Z^T := (Z_1^T , \ldots , Z_r^T) $ a covariable. This model has been studied in \cite{Staedler:10}
The tuning parameter is there taken of order
$\lambda \asymp \sqrt {\log^3 n \log ( p \vee n ) / n } $. The parameters in this model are
$\theta:= ( \pi, \sigma, \beta_1 , \ldots , \beta_r ) $ (and in \cite{Staedler:10}, the penalty is 
$\lambda \| \beta \|_1 = \lambda \sum_{k=1}^r \| \beta_k \|_1$, i.e., it does not include
the parameters $\pi $ and $\sigma$).
The model
is definitely non-linear. 
It is essentially a GLM albeit that there are $r$ linear functions involved
instead of just one, and there are the further parameters $\pi$ and
$\sigma$. We call such a model an {\it extended} GLM.
The contraction inequality will not help us anymore in this case, but as we will see
in Subsection \ref{non-linear.section}, the
generic chaining argument gives
a multivariate version of the contraction theorem. This leads to the reduced choice
$\lambda \asymp \sqrt {\log p/ n } $.

Another situation is where $\rho_{\theta}$ is a general non-linear function.
In that case, we will restrict ourselves to the medium-dimensional
situation with $p $ sufficiently smaller than $n$. Again the generic
chaining bound can be used. 

Our results rely on the following condition.

\begin{condition}\label{component-Lipschitz.condition} (Componentwise Lipschitz condition)
There exist functions $\{ \psi_j \}$  
($\psi_j : {\cal X} \times \{ 1 , \ldots , n \} \rightarrow \R $) and constants $\{ c_{i, \theta}\}$ such that or all $\theta$ and $\tilde \theta $ in $\Theta_*$
$$
\biggl | [\rho_{\theta} (X_i) - c_{i, \theta} ]- [\rho_{\tilde \theta } (X_i) - 
c_{i, \tilde \theta} ]  \biggr | \le \sum_{j=1}^p | \theta_j - \tilde \theta_j | \psi_j (X_i, i) , \ \forall \ i  . $$
\end{condition}

The constants $\{ c_{i, \theta} \}$ will generally be either all zero, or equal to the expectation
$c_{i, \theta} = \EE \rho_{\theta}(X_i) $. 

Generic chaining gives a bound $\gamma_2$ (following the notation of \cite{talagrand2005generic})
for the supremum of stochastic processes (see Theorem \ref{majorizing.theorem}). This bound $\gamma_2$ is defined
by the geometry of the index set of the process. By Sudakov's minoration
$\gamma_2$ is also a lower bound in the case of Gaussian processes. This is the argument
we will use. It means that we need not directly calculate $\gamma_2$ but instead
obtain an upper bound for free. Nevertheless, it would
be of interest to directly bound $\gamma_2$ using geometric arguments (Talagrand's research problem 2.1.9 in
\cite{talagrand2005generic}). 
The Dudley bound (see \cite{dudley1967sizes} or \cite{dudley2010sizes}) results in additional (and hence superfluous) $(\log n )$-factors
(see Section \ref{geometry.section}).

We remark that the bounds are based on arguments for Gaussian processes,
and in fact on the behavior of maxima of i.i.d.\ Gaussians.
This is so to speak the worst case: the bounds are here the largest.
In particular for random variables which are highly dependent,
one may have smaller bounds. Moreover, in the statistical application
of $\ell_1$-regularized estimation, strong dependencies may lead
to choosing the tuning parameter $\lambda$ of much smaller order
than $\sqrt {\log p / n}$. This is explained in \cite{Lederer:11} for the case
of the Lasso. It means that even when the result
$$\sup_{\theta \in \Theta_M (\theta^*) } | (P_n - P) ( \rho_{\theta} - \rho_{\theta^*} ) |
= {\mathcal O}_{\PP}  \biggl ( \sqrt { \log p \over n } \biggr ) M , $$
leaves no room for improvement, there are situations where the
choice $\lambda \asymp \sqrt {\log p / n} $ is much too large.
We will not address this issue here but refer to \cite{Lederer:11} . 

That generic chaining arguments can be used to theoretically show
that $\lambda \asymp \sqrt {\log p / n }$ is appropriate
is perhaps of little practical value. One may argue for example
that cross-validation will rather be used in practice, instead of
a theoretical value. Our finding is 
primarily interesting from a theoretical point of view. 

Generic chaining plays an important role in the statistics literature,
for example to empirical risk minimization (\cite{bartlett2006empirical}), PAC-Bayesian learning
(\cite{audibert2007combining}), and the Lasso with random design (\cite{bartlett2009}).
We believe the application in this paper, addressing the
choice of the tuning parameter $\lambda$ in $\ell_1$-regularization for M-estimators, 
is  an nice opportunity to clearly demonstrate  the elegance of Talagrand's
approach. 

\subsection{Organization of the paper}
In Section \ref{oracle.section}, we review the basic oracle inequality for the
$\ell_1$-penalized M-estimator. This purpose of this section is to highlight the
role of the supremum
$$\sup_{\theta \in \Theta_M(\theta^*)} 
| ( P_n - P) ( \rho_{\theta} - \rho_{\theta^*} ) | . $$
The proofs of Theorems \ref{oracle.theorem} and \ref{oracle2.theorem} follow closely \cite{BvdG2011}, and
are given for completeness in Section \ref{proofs.section}.
In Section \ref{expectation.section} we show that
$$ \EE\biggl (  \biggl [ \sup_{\theta \in \Theta_M (\theta^*) }|  P_n^{\varepsilon} ( \rho_{\theta}^c - \rho_{\theta^*}^c )| \biggr ] \biggr \vert 
{\bf X} \biggr ) 
= {\mathcal O}_{\PP} \biggl (  \sqrt {\log p \over n } \biggr ) M. $$
Here, $P_n^{\varepsilon}$ is the symmetrized measure defined in 
Section \ref{symmcon.section} and ${\bf X} := (X_1 , \ldots , X_n)$. Moreover, $\rho_{\theta}^c (X_i, i ) =
\rho_{\theta} (X_i) -c_{i, \theta}$, with the constants $c_{i, \theta}$ as in Condition \ref{component-Lipschitz.condition}.
Section  \ref{symmcon.section} summarizes why bounds
on the conditional mean of the symmetrized process suffice: they lead to exponential probability
inequalities using a deviation inequality of \cite{Massart:00a}. Section
\ref{geometry.section} gives the details concerning generic chaining and
a consequence concerning the geometry of $\ell_1$-balls. It summarizes some results in
\cite{talagrand2005generic} and makes a comparison with Dudley's entropy bound.

\section{The oracle inequality}\label{oracle.section}

We let for $\theta$ and $\theta^*$ in $\Theta$, 
$$Y(\theta, \theta^*) := (P_n - P) (\rho_{\theta} - \rho_{\theta^*} ) .$$

In this section, we show why bounds for $\sup_{\theta \in \Theta_M (\theta^*)} | Y( \theta , \theta^*) |$
can be used to choose the tuning parameter $\lambda$ and arrive at an oracle inequality for
the $\ell_1$-regularized M-estimator $\hat \beta$. 
The line of reasoning is as in \cite{BvdG2011}.
Define the excess risk
$${\cal E} ( \theta ; \theta_0) := P ( \rho_{\theta} - \rho_{\theta_0} ) . $$

The following condition quantifies the curvature of ${\cal E} ( \theta ; \theta_0)$
around its minimizer $\theta_0$. 

\begin{condition} \label{margin.condition} (Margin condition) 
We say that the margin condition holds for
all $\theta \in \Theta_M (\theta^*) $ if
for some norm $\tau$ on $\Theta$, and some strictly convex
non-negative function $G$, satisfying $G(0)=0$, 
$${\cal E} ( \theta ; \theta^0) \ge G( \tau( \theta-\theta_0) ) , \ \forall \ \theta \in \Theta_M (\theta^*) . $$
\end{condition}

\begin{definition}\label{conjugate.definition} (Convex conjugate) Let $G$ be a strictly convex non-negative
function with $G(0) =0$. The convex conjugate of $G$ is
$$H(v) := \sup_{u \ge 0 } \biggl \{ uv - G(u) \biggr \} , \ v \ge  0 . $$

\end{definition}

For sets $S$ and vectors $\theta \in \R^p$ we let
$$\theta_{j,S} := \theta_j {\rm l} \{ j \in S \} , \ j=1 , \ldots , p . $$

\begin{definition}\label{compatibility.condition} (Effective sparsity)
Let 
$$\delta(L,S) := \min \{ \tau(\theta):\ \| \theta_S\|_1 = 1 , \| \theta_{S^c} \|_1 \le L \}  . $$
Then  $\Gamma^2 (L,S) := 1/ 
\delta^2 (L,S)$ is called the effective sparsity (of the set $S$). 
\end{definition} 

Following \cite{vandeG07}, we call  $\phi^2 (L,S) := |S| \delta^2 (L,S)$ the
compatibility constant (for the set $S$). If it is not too small, 
the norms $\tau$ and the $\ell_1$-norm $\| \cdot \|_1$ are ``compatible"
with each other. 

We define for some constant $\lambda_0$, the set
$${\cal T}_M (\theta^*) := \{ |Y(\theta, \theta^* ) |\le \lambda_0 \| \theta - \theta^*\|_1 \vee \lambda_0^2  , \ \forall \ \theta \in \Theta_M (\theta^*)  \} ,$$
and let ${\cal T} (\theta^* ) := {\cal T}_{\infty} ( \theta^* ) $ and $\Theta_{\infty} (\theta^*) = \Theta $. 

Our task in Sections \ref{symmcon.section} and \ref{expectation.section} is to show that with
$\lambda_0 \asymp \sqrt {\log p / n }$, the set ${\cal T}_M (\theta^*)$
has large probability (for any $\theta^*$ and suitable $M $). 

We first give in Theorem \ref{oracle.theorem} a result where the margin assumption is assumed to hold "globally".
We then refine this in Theorem \ref{oracle2.theorem} to local conditions for the convex case.

\begin{theorem} \label{oracle.theorem} 
Let $\lambda > \lambda_0$. 
Assume Condition \ref{margin.condition} (the margin condition) for all $\theta \in \Theta$.
Let $H$ be the convex conjugate of $G$.
If  $\theta^0 \in \Theta$, we have on ${\cal T} (\theta_0)$, for all $0< \delta < 1 $, 
\begin{equation}\label{bound1}
(1- \delta ) {\cal E} ( \hat \theta ; \theta_0) + (\lambda- \lambda_0) \| \hat \theta - \theta^0  \|_1\le 
\delta H \biggl ( { 2 \lambda \Gamma(L , S_0) \over \delta } \biggr )\vee 2 \lambda^2  ,
 \end{equation}
 with 
 $
 L= { (\lambda+ \lambda_0) / (\lambda - \lambda_0) } $.
Moreover, for all $0 < \delta < 1$ and all $\theta^* \in \Theta$, 
on ${\cal T} (\theta^*)$, 
\begin{equation} \label{bound2}
(1-  \delta ) {\cal E} ( \hat \theta ; \theta_0) + (\lambda- \lambda_0) \| \hat \theta- \theta^*  \|_1 $$ $$\le
2 \delta H \biggl ( { 4 (1+  \delta )  \lambda \Gamma(L_{\delta}, S_*) \over \delta^2  } \biggr )\vee 2 \lambda^2 +
(1+ \delta ) {\cal E} ( \theta^* ; \theta_0) , 
\end{equation}
with $L_{\delta} = 2 \left (  { (1+  \delta) / \delta } \right ) \left ( {  (\lambda + \lambda_0 ) /( \lambda - \lambda_0 )} \right )$.
Here $S_* := \{ j : \ \theta_j^* \not= 0 \} $ is the support set of $\theta^*$. 
\end{theorem} 

The proof of Theorem \ref{oracle.theorem} is given in Section \ref{proofs.section}.
\begin{remark}\label{theta*.remark}
With the above result, one can define the best sparse approximation as a solution $\theta^*$ of the minimization
$$\min_{\theta \in \Theta} \left \{ 2 \delta H \biggl ( { 4 (1+  \delta )  \lambda \Gamma(L_{\delta}, S_{\theta}) \over \delta^2  } \biggr )\vee 2 \lambda^2 +
(1+ \delta ) {\cal E} ( \theta ; \theta_0)  \right \} , $$
where $S_{\theta}:= \{ j : \ \theta_j \not= 0 \}$.
\end{remark}

The next theorem assumes convexity and then needs the margin condition only in a neighborhood of $\theta^*$.

\begin{theorem} \label{oracle2.theorem} Let $\lambda > \lambda_0$.
Let $\theta_*$ be the smallest set containing $\Theta$ and 
suppose that the map
$\theta \mapsto \rho_{\theta} $, $\theta \in \Theta_* $ is convex. 
Let $(\lambda- \lambda_0) M_0$ and $(\lambda- \lambda_0) M_*$ be the bounds given in  the
right hand side of (\ref{bound1}) and (\ref{bound2})
respectively, i.e., 
$$M_0 :=  {\delta  \over \lambda - \lambda_0}  \biggl \{ H \biggl ( { 2 \lambda \Gamma(L , S_0) \over \delta } \biggr )\vee 2 \lambda^2 \biggr \}  ,$$
with 
 $
 L= { (\lambda+ \lambda_0 )/ (\lambda - \lambda_0) } $.
and
$$M_* := { 1  \over \lambda- \lambda_0 } \biggl \{ 2 \delta H \biggl ( { 4(1+  \delta )  \lambda \Gamma(L_{\delta}, S_*) \over  \delta^2  } \biggr )
\vee 2 \lambda^2 +
(1+ \delta ) {\cal E} ( \theta^* ; \theta_0)\biggr \} , $$
with 
$L_{\delta} =2 \left (  {( 1+  \delta )/ \delta } \right ) \left ( {  (\lambda + \lambda_0) /( \lambda - \lambda_0 )} \right )$.
Here, $H$ is a strictly convex increasing function with $H(0)=0$. 

If $\theta^0 \in \Theta$ and the margin condition holds for all
$\theta \in \Theta_{2M_0} (\theta_0)$, with $G$ the convex conjugate of $H$, then again
on ${\cal T}_{2M_0}  (\theta_0)$, 
$$(1- \delta ) {\cal E} ( \hat \theta ; \theta_0) + (\lambda- \lambda_0) \| \hat \theta - \theta^0  \|_1\le 
(\lambda- \lambda_0) M_0  $$
 For general $\theta^*$, if the margin condition holds for all $\theta \in \Theta_{2M_*} (\theta^*)$,
with $G$ the convex conjugate of $H$, then again
 on ${\cal T}_{2 M^*}  (\theta^*)$, 
$$(1- 2 \delta ) {\cal E} ( \hat \theta ; \theta_0) + (\lambda- \lambda_0) \| \hat \theta- \theta^*  \|_1 
\le (\lambda- \lambda_0 ) M_* . $$ \end{theorem}

The proof of Theorem \ref{oracle2.theorem} is also in Section \ref{proofs.section}.

\begin{remark}\label{peeling.remark}
To handle the set ${\cal T}_M (\theta^*)$ we prove in Sections \ref{symmcon.section} and
\ref{expectation.section} that with $\lambda_0 \asymp \sqrt {\log p / n }$, with large probability
$$\sup_{\theta \in \Theta_M (\theta^*) } | Y( \theta , \theta^* ) | \le \lambda_0 M $$
for all $M \le {\rm const.}$ and then apply the peeling device
(the latter being detailed in Subsection \ref{peeling.section}). However, as is clear from the proof of Theorem
\ref{oracle2.theorem}, one can refrain from peeling in the convex case,
because one already places oneself in a suitable neighborhood of $\theta^*$
 (see also
\cite{vandeG07} and \cite{vandeG08}).

\end{remark}

\begin{remark} \label{convex.remark} The models considered in
\cite{Staedler:10} and \cite{schelldorfer2011estimation} are not convex. There, 
the margin condition holds in a bounded neighborhood and these bounds are imposed
on the parameters. Then the peeling device is invoked.
\end{remark}

\section{Symmetrization, contraction and deviation inequalities, and the peeling device}\label{symmcon.section}
We write the sample as ${\bf X} := ( X_1 , \ldots , X_n)$.
Let
$\varepsilon_1 , \cdots , \varepsilon_n$ be a Rademacher sequence independent of $\bf X$. 
For constants $\{ c_{i, \theta} \}$ (which we will choose as in
Condition \ref{component-Lipschitz.condition}), we define 
$$ \rho_{\theta}^c (X_i, i) = \rho_{\theta} (X_i) - c_{i , \theta} , \ i=1 , \ldots , n, $$ and
the symmetrized empirical process
$$P_n^{\varepsilon} \rho_{\theta}^c := {1 \over n} \sum_{i=1}^n[  \rho_{\theta} (X_i) - c_{i, \theta} ] \varepsilon_i , \ \theta \in \Theta_*, $$
and we let 
$$Y^{\varepsilon} ( \theta , \theta^*)  := P_n^{\varepsilon} ( \rho_{\theta}^c - \rho_{\theta^* }^c ) , \ \theta \in \bar \Theta. $$
For a function $g : {\cal X} \times \{ 1 , \ldots n \}$, we use the notation
$$\| g \|_n^2 := {1 \over n} \sum_{i=1}^n g^2 (X_i, i)  , \  \|  g \|^2 := {1 \over n} \sum_{i=1}^n  \EE g^2 (X_i, i) 
 . $$

In this section, we summarize the arguments that show that up to constants, one can reduce the problem
of deriving probability inequalities for the process $Y( \theta , \theta^*)$ to studying the symmetrized process
$Y^{\varepsilon} ( \theta , \theta^* ) $. In fact, we only need bounds for the conditional expectation
$$ E_n := \EE\biggr ( \biggl [  \sup_{\theta \in \Theta_M(\theta^* ) } | Y^{\varepsilon} ( \theta , \theta^* )| \biggr ]  \biggr \vert {\bf X} 
\biggr ) . $$
 Alternatively, one can use
direct arguments in certain regression problems (with
sub-Gaussian errors) or invoking or example Bernstein's inequality (but then one has
to adjust Sudakov's minoration argument to the case of independent Gamma-distributed variables).
We also discuss the peeling device (but as noted in Remark \ref{peeling.remark} this device
is not always needed).

\subsection{Symmetrization} \label{symmetrization.section}

We cite the following result (see \cite{pollard1984convergence}).

\begin{lemma} Let $R := \sup_{\theta \in \Theta_M (\theta^*) } \| \rho_{\theta}^c - \rho_{\theta^* }^c \|  $
and let $t \ge 4 $.
Then
$$\PP \biggl ( \sup_{\theta \in \Theta_M (\theta^*)  } | Y ( \theta , \theta^* ) | > 4R \sqrt {2t\over n}   \biggr ) \le 4
\PP \biggl (  \sup_{\theta \in \Theta_M (\theta^*) } | Y^{\varepsilon}  ( \theta , \theta^* ) | >R  \sqrt {2 t \over n}  \biggr ). $$

\end{lemma} 

\subsection{Contraction} \label{contraction.section} 
Suppose that for all $\theta , \tilde \theta \in \Theta_*$, 
$$\biggl |  \rho_{\theta}^c ( X_i, i)-  \rho_{\tilde \theta}^c ( X_i, i ) \biggr  | \le
| f_{\theta} (X_i, i ) - f_{\tilde \theta } (X_i, i ) | ,  \ \forall \ i ,$$
for some functions
$f_{\theta} : {\cal X} \times \{ 1 , \ldots , n \} \rightarrow \R $, $\theta \in \Theta_*$. 

By the contraction inequality of \cite{Ledoux:91},  
$$E_n := \EE \biggl ( \biggl [ \sup_{\theta \in \Theta_M (\theta^*) } | Y^{\varepsilon}  (\theta , \theta^* ) | \biggr ] \biggr \vert {\bf X} \biggr ) 
\le  2 \EE \biggr ( \biggl [ \sup_{\theta \in \Theta_M (\theta^*) } 
| X^{\varepsilon} (\theta, \theta^* ) |\biggr ]  \biggr \vert {\bf X} \biggr ) , $$
with
$$X^{\varepsilon} (\theta, \theta^* ) := P_n^{\varepsilon} ( f_{\theta} - f_{\theta^*} ) := {1 \over n} 
\sum_{i=1}^n \varepsilon_i ( f_{\theta} (X_i, i) - f_{\tilde \theta} (X_i , i) ) . $$

\subsection{A deviation inequality} \label{deviation.section}
Write
$$R_n  := \sup_{\theta \in \Theta_M (\theta^*) } \| \rho_{\theta}^c - \rho_{\theta^*}^c \|_n . $$
We have for all $t > 0$ (see \cite{Massart:00a}), 
$$\PP \biggl ( \biggl [ \sup_{\theta \in \Theta_M (\theta^*) } |Y^{\varepsilon}  (\theta , \theta^*)  |\biggr ]  \ge E_n + R_n \sqrt { 2t \over n } \biggr ) 
\le \exp [-t] . $$
Combining this with the symmetrization result of Section \ref{symmetrization.section}, we obtain the following 
corollary.

\begin{corollary} Let for some $\bar R$ 
$$\sup_{\theta \in \Theta_M (\theta^*) } \| \rho_{\theta}^c - \rho_{\theta^*}^c \| \le \bar R ,$$ 
and let $t \ge 4$. 
 Then for any $\bar E$, 
$$\PP \biggl ( \biggl [ \sup_{\theta \in \Theta_M} |Y (\theta , \theta^*)  |\biggr ]  \ge 8 \bar E   + 4 \bar R \sqrt { 2t \over n } \biggr ) $$ $$
\le 4 \exp [-t] + 4 \PP ( R_n > \bar R  \vee  E_n > \bar E)  .  $$
\end{corollary} 

As for the random variables $E_n$ and $R_n$,  in our context we use Condition \ref{component-Lipschitz.condition}.
Consider first $R_n$.  Condition \ref{component-Lipschitz.condition} yields  by the triangle inequality
$$\sup_{\theta \in \Theta_M (\theta^*) } \| \rho_{\theta}^c - \rho_{\theta^* }^c \|_n \le  M K_n , $$
where
$$ K_n := \max_{1 \le j \le p } \| \psi_j \|_n . $$
Thus, on the set
\begin{equation}\label{T0.equation}
{\cal T}_0 := \biggl \{  \max_{1 \le j \le p } \| \psi_j \|_n \le \bar K \biggr \} , 
\end{equation}
(where $\bar K$ is some constant)
we can bound  the random radii $R_n$ by $ M \bar K$. 
We will see in Section \ref{expectation.section} that a bound for the conditional expectation  
$E_n$ also only involves $K_n$:
$$E_n \le \lambda_0 M K_n , $$
for some constant $\lambda_0 \asymp \sqrt {\log p / n }$. 

In some cases (regression with fixed design) $K_n$ is not random, and the assumption
$\max_{1 \le j \le p } \| \psi_j \|_n \le \bar K $ is a matter of normalization. In other situations,
one can for example apply Bernstein's inequality (\cite{Bennet:62}):

\begin{lemma} \label{Bernstein-psi.lemma} Suppose that the $\psi_j (X_i)$ are uniformly
sub-Gaussian, that is, for some positive constants $L$ and $\tau$ , it holds for all $j$,
$${2 L^2  \over n} \sum_{i=1}^n \biggl [ \EE \exp[ \psi_j^2 (X_i) / L^2 ] -1 \biggr ]  \le \tau^2 . $$
Then for all $t > 0 $, 
$$\PP \biggl ( \max_{1 \le j \le p } \biggl | \| \psi_j\|_n^2 -  \| \psi_j \|^2 \biggr | \ge{ 2 \tau L } \sqrt {2(t + \log p ) \over  n} + { 2L(t + \log p)   \over n}  \biggr ) \le 2 \exp[-t] . 
$$

\end{lemma}

{\bf Proof.} The sub-Gaussianity implies that for all $m \in \{ 1, 2,3, \ldots \} $, 
$${1 \over n} \sum_{i=1}^n \EE | \psi_j^2 (X_i) |^m/ n  \le{ L^{2m} m!  \over n} \sum_{i=1}^n \biggl [ \EE \exp[ \psi_j^2 (X_i) / L^2 ] -1\biggr ]    \le {m! \over 2}L^{2(m-1) } \tau^2 .$$
But then 
$${1 \over n} \sum_{i=1}^n \EE | \psi_j^2 (X_i) - \EE \psi_j^2 (X_i) |^m \le { m! \over 2} 2^m L^{2(m-1)} \tau^2 . $$
By Bernstein's inequality (\cite{Bennet:62}), for all $t >0$,
$$ \PP \biggl ( | (P_n - P) \psi_j^2 | \ge { 2 \tau L } \sqrt {2t  \over  n} + { 2Lt  \over n}  \biggr ) \le 2 \exp[-t] ,  $$
and hence, by the union bound, for all $t > 0$, 
$$\PP \biggl ( \max_{1 \le j \le p } | (P_n - P) \psi_j^2 | \ge { 2 \tau L } \sqrt {2(t + \log p ) \over  n} + { 2L(t + \log p)   \over n}  \biggr ) \le 2 \exp[-t] .  $$
\hfill $\sqcup \mkern -12mu \sqcap$

The assumption of sub-Gaussianity is not a necessary condition. One may replace it by an
$m$-th order moment condition, with $p^{2} / n^m $ sufficiently small. However, we then will no longer have exponential
probability inequalities.

\subsection{The peeling device} \label{peeling.section} 
The peeling device goes back to \cite{Alexander:85}, the terminology being introduced in \cite{vandeGeer:00}.
In the present context we can use it in the following form. 

We show in the next section that under certain conditions 
\begin{equation}\label{mean.equation}
E_n \le \lambda_0 M K_n  , 
\end{equation}
where $\lambda_0 \asymp \sqrt {\log p / n }$, and
$K_n := \max_{1 \le j \le p } \| \psi_j \|_n $. 
Then, under Condition \ref{component-Lipschitz.condition} and the sub-Gaussianity assumption of Lemma \ref{Bernstein-psi.lemma}, 
we have for 
all $M $, and  all $t >0$, 
\begin{equation}\label{fixed-M.equation}
\PP \biggl (  \sup_{\theta \in \Theta_M (\theta^*) } | Y( \theta , \theta^* ) | \ge { \lambda_* M \over {\rm e} }  \biggl (1 + K_* \biggl [ \sqrt {t \over \log p } + 
{t \over n} \biggr ] \biggr )  \biggr ) \le
6\exp [-t] . 
\end{equation}
with
$K_*$ depending on $L$ and $\tau$ but not on $n$ and $p$, and $\lambda_* \asymp \sqrt {\log p / n} $.
This follows from (\ref{mean.equation}) and from Subsection \ref{deviation.section}. The constant 6 in the right hand side of inequality (\ref{fixed-M.equation})
comes from a 4 from the symmetrization plus a 2 from Lemma \ref{Bernstein-psi.lemma} (we actually may replace 2 by 1 here
because we only need a one-sided version). 

Once (\ref{fixed-M.equation}) is established, we can invoke the peeling device as follows. Let $\bar M$ be fixed, and let
$M_j := {\rm e}^{-j} \bar M $, $j=0 , \ldots , {p} $. Then for all $t >0$, 
$$\PP \biggl ( \sup_{\theta \in  \Theta_{\bar M} (\theta^*) } { |Y( \theta , \theta^* ) | 
\over  \| \theta - \theta^* \|_1 \vee {\rm e}^{-(p-1)} \bar M } \ge \lambda_* 
\biggl ( 1 + K_*\biggl [ \sqrt {t+ \log p   \over  \log p } + { t  + \log p \over n } \biggr ]  \biggr ) \biggr ) $$
$$\le\sum_{j=1}^p  \PP\biggr  (\sup_{\theta \in \Theta_{M_{j-1} }(\theta^*) } | Y( \theta , \theta^* ) | > \lambda_* M_{j}  \biggl ( 1+  K_*
\biggr [ \sqrt {t+\log p  \over  \log p  } + {t + \log p \over n} \biggr ] \biggr )  \biggr ) $$ $$ \le 6 \exp[ \log p -(\log p   +t) ] \le 6 \exp[-t].  $$

\section{Bounds for the symmetrized process} \label{expectation.section}
In the previous section we argued that the main task is to establish bounds for the expectation of the symmetrized
process, i.e., for
$$\EE\biggl ( \biggl [  \sup_{\theta \in \Theta_M(\theta^*) } | Y^{\varepsilon} ( \theta , \theta^* ) | \biggr ] \biggr \vert {\bf X} \biggr ) . $$
We are looking for bounds of the type (\ref{mean.equation}).  One
can then derive deviation inequalities as shown in Section \ref{symmcon.section}, and hence (as shown in Theorems \ref{oracle.theorem}
and \ref{oracle2.theorem}) theoretical bounds for the
tuning parameter of the $\ell_1$-regularized M-estimator.

\subsection{Linear functions} \label{linear.section}
Lets us briefly recall the linear case.  Let $[\rho_{\theta}(X_i) -c_{i, \theta} ] - [ \rho_{\tilde \theta}(X_i) -c_{i, \tilde \theta}] $ be linear:
$$[\rho_{\theta}(X_i) -c_{i, \theta} ] - [ \rho_{\tilde \theta}(X_i) -c_{i, \tilde \theta}] = \sum_{j=1}^p ( \theta_j - \tilde \theta_j ) \psi_j (X_i , i ), \ i=1 , \ldots , n .  $$
One then clearly has
$$  \EE\biggl ( \biggl [  \sup_{\theta \in \Theta_M(\theta^*) } | Y^{\varepsilon} ( \theta , \theta^* ) | \biggr ] \biggr \vert {\bf X} \biggr )\le M  \| \varepsilon^T \psi/n  \|_{\infty}  .$$
Moreover, by Hoeffding's inequality (\cite{Hoeffding:63}, see also Lemma 14.14 in \cite{BvdG2011})
$$\EE \biggl ( \| \varepsilon^T \psi/n  \|_{\infty} \biggr \vert {\bf X}  \biggr ) \le \sqrt { 2 \log (2p ) \over n }K_n , $$
where $K_n  := \max_{1 \le j \le p } \| \psi_j \|_n $.

\subsection{Generalized linear functions} \label{GLM.section}
Suppose that for all $\theta , \tilde \theta \in \Theta_*$, 
$$\biggl | [ \rho_{\theta} ( X_i) - \EE \rho_{\theta} (X_i) ] - [ \rho_{\tilde \theta} ( X_i) - \EE \rho_{\tilde \theta} (X_i) ]\biggr  | \le
| f_{\theta} (X_i, i) - f_{\tilde \theta } (X_i, i ) | , \ \forall \ i , $$
where  
$f_{\theta} (X_i, i)= \sum_{j=1}^p \theta_j \psi_j (X_i , i)  $, $\theta \in \Theta_*$. 
Then by the contraction inequality of Subsection \ref{contraction.section}, and the arguments
of Subsection \ref{linear.section} for the linear case
$$\EE  \biggl ( \biggl [  \sup_{\theta \in \Theta_M (\theta^* ) }| Y^{\varepsilon}  (\theta, \theta^* ) | \biggr \vert
{\bf X} \biggr ) \biggr ) \le 2M \sqrt {2 \log (2p) \over n } K_n , $$
with $K_n:= \max_{1 \le j \le p } \| \psi_j \|_n $.

\subsection{Extended generalized linear functions} \label{extended-GLM.section} 

\begin{condition} \label{Lipschitz.condition}(Extended GLM condition)
The exist non-negative functions $\{ \psi_{j,k} : j=1 , \ldots p_k, \ k=1 , \ldots , r\}  $ (with $\sum_{k=1}^r p_k = p$) such that for all $\theta $ and $\tilde \theta$ in $\Theta_*$, 
it holds that
$$| [\rho_{\theta}(X_i) - c_{i, \theta} ]  - [ \rho_{\tilde \theta} (X_i) -
c_{i, \tilde \theta} ] | \le \sum_{k=1}^r   | \sum_{j=1}^{p_k}  (\theta_{j,k} - \tilde \theta_{j,k} ) \psi_{j,k} (X_i, i)  | , \ i=1 , \ldots , n .  $$

\end{condition}

\begin{theorem}\label{multidim-contraction.theorem}(Multivariate contraction theorem)
Assume Condition \ref{Lipschitz.condition}.
Let $\xi_{1,k} , \ldots , \xi_{n,k} $, $k=1 , \ldots , r$, be independent ${\cal N} (0, 1)$-distributed random
variables, independent of $X_1 , \ldots , X_n$. Let
$$X_k (\theta, \theta^* )  := {1 \over \sqrt n} \sum_{i=1}^n \sum_{j=1}^{p_k}  ( \theta_{j,k} - \theta_{j,k}^*)  \psi_{j,k}  (X_i,i) \xi_{i,k} ,$$
and
$$X (\theta, \theta^*) := \sum_{k=1}^r X_k (\theta, \theta^* ) = {1 \over \sqrt n} \sum_{i=1}^n \sum_{k=1}^r  \sum_{j=1}^{p_k} ( \theta_{j,k}- \theta_{j,k}^*)  \psi_{j,k} (X_i, i) \xi_{i,k} .  $$
 Then for a universal constant $C$, 
$$\EE \left ( \biggl [ \sup_{\theta \in \Theta_M(\theta^*)  } |Y^{\varepsilon} (\theta, \theta^* ) | \biggr ]  \biggl \vert {\bf X} \right ) \le C 2^{r-1}  \EE\biggr ( \biggl [ \sup_{\theta \in \Theta_M (\theta^*) } X(\theta , \theta^* ) \biggr ]  \biggr \vert {\bf X} \biggr )  . $$
\end{theorem}

{\bf Proof.} We apply Theorem 2.1.1 in \cite{talagrand2005generic}, cited in the present paper as
Theorem \ref{majorizing.theorem}. Note first that
$$\EE \biggl ( |X (\theta , \theta^*) - X (\tilde \theta , \theta^*) |^2 \vert {\bf X} \biggr ) =\sum_{k=1}^r  \| \sum_{j=1}^{p_k}(\theta_{j,k} - \tilde \theta_{j,k} ) \psi_{j,k}\|_n^2 
  .$$
For all $\theta$ and $\tilde \theta$ we have
$$\| \rho_{\theta }^c - \rho_{\tilde \theta}^c \|_n^2 \le  \| \sum_{k=1}^r  | \sum_{j=1}^{p_k} 
( \theta_{j,k} - \tilde \theta_{j,k}) \psi_{j,k} |  \|_n^2 \le 2^{r-1}  \sum_{k=1}^r \| \sum_{j} \theta_{j,k} \psi_{j,k} \|_n^2 .  $$
By Hoeffding's inequality (\cite{Hoeffding:63})
$$\PP \biggl ( | Y^{\varepsilon } (\theta, \theta^*)  - Y^{\varepsilon}(\tilde \theta, \theta^* ) |  \ge \| \rho_{\theta}^c - \rho_{\tilde \theta}^c \|_n  \sqrt {2t} \ \biggr \vert {\bf X} \biggr ) \le 2 \exp[-t] . $$
Hence, using Theorem 2.1.5 in Talagrand's book (\cite{talagrand2005generic}) (see Section \ref{geometry.section},
Theorem \ref{consequence.theorem}), 
we get for a universal constant $C$,
$$\EE \left ( \biggl [ \sup_{\theta \in \Theta_M(\theta^*)  } |Y^{\varepsilon} (\theta, \theta^* ) | \biggr ]  \biggl \vert {\bf X} \right )  \le C 2^{r-1}  \EE\biggl ( \biggl [ \sup_{\theta \in \Theta_M ( \theta^*) } X(\theta , \theta^* ) \biggr ]  \biggr \vert {\bf X} 
\biggr ) . $$
\hfill $\sqcup \mkern -12mu \sqcap$

As a direct consequence (i.e., by bounding the right hand side in Theorem \ref{multidim-contraction.theorem}), we obtain
the bounds of interest for our problem.

\begin{theorem} \label{extended-GLM.theorem} Assume Condition \ref{Lipschitz.condition} and let
$K_n := \max_{j,k} \| \psi_{j,k} \|_n  .$
We have for a universal constant $C$,
$$\EE \left ( \biggl [ \sup_{\theta \in \Theta_M(\theta^*)  } |Y^{\varepsilon} (\theta, \theta^* ) | \biggr ]  \biggl \vert {\bf X} \right )   \le  C2^{r-1}  \sqrt { 2 \log (2p ) \over n} K_n . $$

\end{theorem}

{\bf Proof.} Let $X( \theta , \theta^*)$ be defined as in Theorem 
\ref{multidim-contraction.theorem}. 
As in Subsection \ref{linear.section}, but now for Gaussians instead of a Rademacher sequence,
conditionally on ${\bf X} := (X_1 , \ldots , X_n)$, we have 
$$\EE  \biggl ( \biggl [ \sup_{\theta \in \Theta_M ( \theta^*)  } X (\theta, \theta^* ) \biggr ] \biggl \vert {\bf X} \biggr ) \le M  \sqrt {2 \log(2 p) \over n} K_n. $$
\hfill $\sqcup \mkern -12mu \sqcap$

\subsection{Non-linear functions} \label{non-linear.section}
We now consider the case where the loss $\rho_{\theta}$ is possibly not extended GLM, that is, its dependence on $\theta$ is
strictly non-linear. However, we do assume that it is component-wise Lipschitz in $\theta$, i.e., that Condition \ref{component-Lipschitz.condition}
holds.

Define for $\psi= ( \psi_1 , \ldots , \psi_p)^T $, 
$$\Sigma_n := {1 \over n } \sum_{i=1}^n \psi (X_i, i) \psi^T (X_i, i ) . $$
Let $\underline \Lambda_n^2 $ be the smallest eigenvalue of $\Sigma_n$ and $\bar \Lambda_n^2$ be its largest eigenvalue.
We assume that $\underline \Lambda_n >0$, thus excluding the case $p > n $. 

\begin{theorem} Assume Condition \ref{component-Lipschitz.condition}. For a universal constant $C$, it holds that
$$\EE \left ( \biggl [ \sup_{\theta \in \Theta_M(\theta^*)  } |Y^{\varepsilon} (\theta, \theta^* ) | \biggr ]  \biggl \vert {\bf X} \right )  \  \le  CM
 \sqrt { 2 \log (2p) \over n }   \biggl (  {\bar \Lambda_n / 
\underline \Lambda_n } \biggr ) . $$

\end{theorem} 

{\bf Proof.} Use that
$$\| \sum_{j=1}^k \theta_j \psi_j \|_n^2 \ge  \underline  \Lambda_n^2 \| \theta \|_2^2 , $$
and
$$ \sum_{j=1}^k |\theta_j | \psi_j \|_n^2 \le \bar \Lambda_n^2 \| |\theta | \|_2^2 = \| \theta \|_2^2 . $$
Then apply the same arguments as in Theorem \ref{extended-GLM.theorem}. 
\hfill $\sqcup \mkern -12mu \sqcap$

\section{The geometry of $\ell_1$-balls} \label{geometry.section}
We first describe here the generic chaining bound, specialized to our context and with a notation adjusted to our
setting. Let
$\xi_1 , \ldots , \xi_n $ be independent ${\cal N} (0,1)$-distributed random variables
and ${\cal V} $ be a subset of $\R^n$. Define
$$X_{v} := {1 \over n} \sum_{i=1}^n v_i \xi_i , \ v \in {\cal V} . $$
Moreover, write
$$\| v \|_n^2 := {1 \over n}  \sum_{i=1}^n v_i^2 , \ v \in \R^n . $$
Talagrand (\cite{talagrand2005generic}, Definition 1.2.3) calls a sequence of partitions $\{ {\cal A}_s \}_{s=0}^{\infty}$
of ${\cal V}$ 
{\it admissible}  if it is an increasing sequence (i.e., ${\cal A}_{s+1}$ contains ${\cal A}_s$ for all $s \ge 1$), and
$|{\cal A}_s | \le 2^{2^s} $ for all $s$.  He defines for each $v \in {\cal V}$ and each $s$,  the set 
$A_s (v)$ as the unique element of ${\cal A}_s$ that contains $v$, and $\Delta (A_s (v))$ as the diameter of
$A_s (v)$. He writes
$$\gamma_2 ({\cal V} , \| \cdot \|_n):= \inf \sup_{v \in {\cal V}} \sum_{s \ge 0 } 2^{s/2} \Delta ({\cal A}_s (v))  ,$$
where the infimum is taken over all admissible partitions.

\begin{theorem} (The majorizing measure theorem, see \cite{talagrand2005generic}, Theorem 2.1.1)
\label{majorizing.theorem}
For some universal constant $ C$, we have
$${1 \over C} \gamma_2 ( {\cal V} , \| \cdot \|_n ) \le \EE \biggl [ \sup_{v \in {\cal V} } X_v \biggr ] \le  C
\gamma_2 ( {\cal V} , \| \cdot \|_n ) . $$
\end{theorem}

Talagrand derives the lower bound in the above theorem from Sudakov's minoration argument.
As a consequence, Talagrand presents the following result.

\begin{theorem} (\cite{talagrand2005generic}, Theorem 2.1.5)\label{consequence.theorem}
Let $\{ Y_v: \ v \in {\cal V} \}$ be a stochastic process that satisfies for all $t>0$
$$\PP \biggl ( | Y_v - Y_{\tilde v} | \ge \sqrt t \biggr ) \le 2 \exp \biggl [ - { t \over \| v-\tilde v \|_n^2 } \biggr ] , \ \forall \ v , \tilde v  \in {\cal V}. $$
Then for a universal constant $C$, we have
$$\EE \biggl [ \sup_{v , \tilde v \in {\cal V}} | Y_v - Y_{\tilde v} | \biggr ] \le C \EE\biggl [ \sup_{v , \tilde v \in {\cal V}} | X_v - X_{\tilde v} | \biggr ] . $$
\end{theorem}

Let us compare here the situation with Dudley's entropy bound (\cite{dudley1967sizes}). We formulate
it using chaining along a tree, as in \cite{BvdG2011}, Subsection 14.12.4, or \cite{vandeGeer:11}. Define
$R_n:= \sup_{v \in {\cal V} } \| v \|_n $.
Let for each $s \in \{ 0, 1, \ldots , S\}$,
$\{ v_j^s \}_{j=1}^{N_s} \subset {\cal V} $ be a minimal $2^{-s}R_n$-covering set of ${\cal V}$, that is, for all $v \in {\cal V}$
and all $s$
there is a $v_j^s$ such that $\| v - v_j^s \|_n \le 2^{-s} R_n$. Then for all $v$, we can find a
end node $v^{S}  \in \{ v_j^{S} \}$ such that $\| v- v^{S}  \|_n \le 2^{S}R_n$,
and for each end node $v^{S} \in \{ v_j^{S} \} $ one can find a branch
$\{ v^0 , \ldots , v^S \}$ such that $\| v^{s} - v^{s-1} \|_n \le 2^{-{s-1}}R_n$ for all
$s=1, \ldots , S$. 
Moreover, we can write (with $X_{v^0} =0 $)
$$X_v = \sum_{s=0}^{S} ( X_{v^{s}} - X_{v^{s-1}}) + X_v - X_{v^S} . $$
Invoking
$$| X_v - X_{v^S} | \le 2^{-S} R_n \sqrt {\sum_{i=1}^n \xi_i^2 } / n , $$
one arrives at Dudley's bound
\begin{equation}\label{Dudley.equation}
\EE \biggl [ \sup_{v \in {\cal V}} X_v \biggr ] \le \sum_{s=0}^S 2^{-(s-1)} R_n \sqrt { 2 \log (2 N_s ) \over n } +
{ 2^{-S} R_n  } . 
\end{equation}

Consider now a special case. We let $ \{ \psi_j \}_{j=1}^p $ be $p$ vectors
in $\R^n$, and let
$${\cal V} := \{ \sum_{j=1}^p \theta_j \psi_j : \ \| \theta \|_1 \le 1 \}  .$$
Let $K_n := \max_{1 \le j \le p } \| \psi_j \|_n $.

The following lemma rephrases the first part of Theorem 2.1.6 in \cite{talagrand2005generic}.
We present a short proof to show that it is again based on the dual norm inequality (\ref{dual-norm.equation}).
\begin{lemma} \label{gamma2.lemma} It holds for some universal constant $C$ that
$$\gamma_2 ( {\cal V}, \| \cdot \|_n ) \le C \sqrt { 2 \log (2p) \over n } K_n . $$
\end{lemma}

{\bf Indirect Proof.} 
Clearly, by the dual norm inequality
$$ \sup_{v \in {\cal V} } X_v = \sup_{\| \theta \|_1 \le 1 } 
{1 \over n} \sum_{i=1}^n \sum_{j=1}^p \theta_j \psi_{i,j} \xi_i =\max_{1 \le j \le p }   \biggl | {1 \over n} 
 \sum_{i=1}^n \psi_{i,j} \xi_i \biggr | . $$
 Hence,
 $$ \EE \biggl [ \sup_{v \in {\cal V} } X_v\biggr ] \le \EE \max_{1 \le j \le p }    \biggl | {1 \over n} 
 \sum_{i=1}^n \psi_{i,j} \xi_i \biggr | \le \sqrt { 2 \log (2p ) \over n } K_n . $$
 The result now follows from Theorem \ref{majorizing.theorem}.
\hfill $\sqcup \mkern -12mu \sqcap$

In his book, Talagrand now poses the research question to prove Lemma \ref{gamma2.lemma}
directly (\cite{talagrand2005generic}, Research problem 2.1.9).
We claim that this cannot be done by applying Dudley's bound. Our reasoning is as follows.
Using Theorem 6.2 in \cite{pollard1990empirical} (see also \cite {vanderVaart:96}, Lemma 2.6.11, or  \cite{BvdG2011}, Lemma 14.28),
we see that
\begin{equation}\label{entropy.equation}
\log (2N_s)  \le 2^{2s} \log (4p) , \ \forall \ s . 
\end{equation}
Insert this in (\ref{Dudley.equation}) with the bound $R_n \le K_n$, to find that
$$\EE \biggl [ \sup_{v \in {\cal V}} X_v \biggr ] \le 2 (S+1) K_n \sqrt {2 \log (4p) \over n} + 2^{-S} K_n . $$
Minimizing this over $S$ gives a bound of order $(\log n) \sqrt {\log p / n}  K_n $. In other words
(assuming the entropy bound (\ref{entropy.equation}) is up to constants tight, which we believe it is)
invoking Dudley's bound instead of generic chaining leads to a
redundant $(\log n)$-factor. Apparently, Dudley's bound does not fully capture
the geometry of $\ell_1$-balls.

\section{Concluding remarks}
This paper combines results in literature concerning symmetrization,
contraction, deviation inequalities and chaining. Their application in statistical
theory has been highlighted by \cite{Massart:00b}. We have added now a new application, where
generic chaining allows one to remove additional $\log n$ factors. For example, we
have improved the choice $\lambda \asymp \sqrt {\log^3 n \log ( p \vee n ) / n }$ in \cite{Staedler:10}
to $\lambda \asymp \sqrt {\log p / n }$.  The geometric arguments to bound $\gamma_2$ in the case
of convex hulls are still to be developed. Somehow, the generic chaining bound
$\gamma_2$ better exploits  the impossibility to play cat and mouse.

\section{Proofs of Theorems \ref{oracle.theorem} and \ref{oracle2.theorem}} \label{proofs.section}
{\bf Proof of Theorem \ref{oracle.theorem}.} 
The Basic Inequality says that
$${\cal E} ( \hat \theta ; \theta_0) + \lambda \| \hat \theta \|_1 \le
Y(\hat \theta, \theta^*) + \lambda \| \theta^* \|_1   + {\cal E} ( \theta^* ; \theta_0). $$
Hence on ${\cal T}(\theta^*)$,
$$ {\cal E} ( \hat \theta ; \theta_0) + \lambda \| \hat \theta \|_1 \le \lambda_0 \| \hat \theta - \theta^* \|_1\vee \lambda_0^2  +
\lambda \| \theta^* \|_1 + {\cal E} ( \theta^* ; \theta_0)  . $$

If $\| \hat \theta - \theta^* \|_1 \le \lambda_0$, we get
$${\cal E} ( \hat \theta ; \theta^0) + (\lambda- \lambda_0) \| \hat \theta - \theta^* \|_1 \le 
2 \lambda^2 + {\cal E} ( \theta^* ; \theta^0) . $$
Hence in the rest of the proof, we can assume $\| \hat \theta - \theta^* \|_1 \ge \lambda_0$.

For $\theta^*= \theta^0$, we get
$$ {\cal E} ( \hat \theta ; \theta_0) + ( \lambda - \lambda_0)  \| \hat \theta_{S_0^c}  \|_1 \le 
(\lambda + \lambda_0 ) \| \hat \theta_{S_0}  - \theta^0 \|_1, $$
which gives for any $0 < \delta < 1$, 
$$ {\cal E} ( \hat \theta ; \theta_0) + ( \lambda - \lambda_0)  \| \hat \theta - \theta^0   \|_1 \le 
2 \lambda \| \hat \theta_{s_0}  - \theta^0 \|_1$$
$$ \le 2 \lambda \Gamma (L, S_0) \tau ( \hat \theta - \theta^0 ) $$
$$ \le \delta {\cal E} ( \hat \theta ; \theta^0 ) + \delta H \biggl (
{ 2 \lambda \Gamma ( L , S_0 ) \over \delta } \biggr ) . $$

For general $\theta^*$, we get
$${\cal E} ( \hat \theta ; \theta_0) + ( \lambda- \lambda_0) \| \hat \theta_{S_*^c} \|_1 \le
( \lambda+\lambda_0) \| \hat \theta_{S_*}  - \theta^* \|_1 + {\cal E} ( \theta^* ; \theta_0) . $$
If $( \lambda+ \lambda_0) \| \hat \theta_{S_*} - \theta^* \|_1 \le \delta {\cal E} ( \theta^* ; \theta_0) $,
we obtain
$${\cal E} ( \hat \theta ; \theta_0) + (\lambda - \lambda_0)  \| \hat \theta_{S_*^c}  \|_1 \le (1+ \delta )  {\cal E} ( \theta^* ; \theta_0). $$
And then, using $\lambda - \lambda_0 \le \lambda + \lambda_0$,
$${\cal E} ( \hat \theta ; \theta_0) + (\lambda - \lambda_0)  \| \hat \theta - \theta^*  \|_1 \le (1+ 2 \delta )  {\cal E} ( \theta^* ; \theta_0). $$
If $( \lambda+ \lambda_0) \| \hat \theta_{S_*} - \theta^* \|_1 \ge \delta {\cal E} ( \theta^* ; \theta_0) $,
we obtain
$${\cal E} ( \hat \theta ; \theta_0) + ( \lambda- \lambda_0) \| \hat \theta_{S_*^c} \|_1 \le
{ 1+ \delta \over \delta } ( \lambda + \lambda_0) \| \hat \theta_{S_*}  - \theta^* \|_1 ,$$
and hence
$${\cal E} ( \hat \theta ; \theta_0) + ( \lambda- \lambda_0) \| \hat \theta - \theta^* \|_1 \le
{ 1+ 2 \delta \over \delta } ( \lambda + \lambda_0) \| \hat \theta_{S_*}  - \theta^* \|_1 ,$$
$$\le { 1+ 2 \delta \over \delta }
(\lambda + \lambda_0) \Gamma (L_{\delta} S_*) \tau (\hat \theta- \theta^* )  + {\cal E} ( \theta^* ; \theta^0) $$
$$ \le 4 \delta H \biggl ( { (1+ 2 \delta ) ( \lambda+ \lambda_0) \Gamma(L_{\delta} , S_*) \over 2 \delta^2  } \biggr ) +
\delta {\cal E} ( \hat \theta ; \theta^0) + (1+ \delta)  {\cal E} ( \theta^* ; \theta^0) . $$
It follows hat
$$(1- 2 \delta ) {\cal E} ( \hat \theta ; \theta_0) + (\lambda- \lambda_0) \| \hat \theta- \theta^*  \|_1 \le
4 \delta H \biggl ( { (1+ 2 \delta ) ( \lambda+ \lambda_0) \Gamma(L_{\delta}, S_*) \over 2 \delta^2  } \biggr )$$ $$+
(1+ 2\delta ) {\cal E} ( \theta^* ; \theta_0) . $$
Finally simplify the expression using $\lambda + \lambda_0 \le 2 \lambda$, and replacing $2\delta$ by $\delta$.

\hfill$ \sqcup \mkern -12mu \sqcap$

{\bf Proof of Theorem \ref{oracle2.theorem}.} We only describe the case $\theta^* = \theta^0$, the case $\theta^* \not= \theta^0$ following by the
same arguments.
Repeat the proof of Theorem \ref{oracle.theorem} with $\hat \theta$ replaced by
$\tilde \theta := t \hat \theta + (1- t) \theta^0 $, where
$$t: = { 2M_0 \over 2M_0 + \| \hat \theta - \theta^0 \|_1 } . $$
Note that $\| \tilde \theta - \theta^0 \|_1 \le 2 M_0$. By the proof of Theorem \ref{oracle.theorem},
we obtain that actually $\| \tilde \theta- \theta^0 \|_1 \le M_0$ on ${\cal T}_{2M_0} (\theta_0)$.
But this implies $\| \hat \theta - \theta_0 \|_1 \le 2 M_0$. Now, repeat the proof again, knowing
that $\| \hat \theta - \theta^0\|_1 \le 2 M_0 $.

\hfill$ \sqcup \mkern -12mu \sqcap$





\bibliographystyle{plainnat}
\bibliography{reference}








\end{document}